# GLOBAL IRREGULARITY
# FOR MILDLY DEGENERATE ELLIPTIC OPERATORS

MICHAEL CHRIST

ABSTRACT. Examples are given of degenerate elliptic operators on smooth, compact manifolds that are not globally regular in $C^\infty$. These operators degenerate only in a rather mild fashion. Certain weak regularity results are proved, and an interpretation of global irregularity in terms of the associated heat semigroup is given.

## 1. INTRODUCTION

Let $L$ be a linear partial differential operator, defined on a compact $C^\infty$ manifold $M$ without boundary.

**Definition.** $L$ is said to be globally regular in $C^\infty$ if
$$\big[u \in \mathcal{D}'(M) \text{ and } Lu \in C^\infty(M)\big] \implies u \in C^\infty(M).$$

Global $C^\infty$ regularity should be contrasted with $C^\infty$ hypoellipticity, which means that for any open set $V \subset M$, we have $u \in C^\infty(V)$ for any $u \in \mathcal{D}'(V)$ such that $Lu \in C^\infty(V)$. The former is implied by the latter, but is in general a far weaker property. As an example, consider $L = \partial_{x_1} + \alpha \partial_{x_2}$ on $M = \mathbb{T}^2$, the two-dimensional torus, where $\alpha$ is any real constant. Such a real vector field is never hypoelliptic, yet it is globally regular in $C^\infty$ for almost every $\alpha$ (in the sense of Lebesgue measure), as may be seen by comparing Fourier coefficients of $u$ and $Lu$. These examples are a bit artificial, but there exist other, more natural examples of operators globally regular in $C^\infty$, yet not hypoelliptic. One class of such examples will be described in the first remark of Section 9.

Consider the case where $L$ is a second-order operator, with $C^\infty$ coefficients and nonnegative principal symbol, that is selfadjoint on $L^2(M, \sigma)$, where the measure $\sigma$ has a $C^\infty$, nowhere vanishing density with respect to Lebesgue measure in any local coordinate system. In order to exclude various degenerate examples in which global regularity fails to hold for trivial reasons, we assume that:

*Date*: November 29, 1995.
Research supported by NSF grant DMS-9306833 and at MSRI by NSF grant DMS-9022140.

(1.1) $L$ is elliptic at every point of an open dense subset of $M$;
(1.2) there exists $C < \infty$ such that $\|u\| \leq C\|Lu\|$ for all $u \in C^2(M)$, and
(1.3) for any points $p, q \in M$ there exists a subunit path joining $p$ to $q$.

By a subunit path we mean a Lipschitz continuous function $\gamma : [0, R] \to M$, for some $R < \infty$, such that, for almost every $t$, the tangent vector $\dot\gamma(t)$ has length less than or equal to 1 with respect to the degenerate Riemannian metric associated to $L$. Such a path is said to join $p$ to $q$ if $\gamma(0) = p$ and $\gamma(1) = q$. In local coordinates $x$, if $L = -\sum a_{ij} \partial^2_{x_i x_j}$ plus lower-order terms and if $\gamma(t) = (\gamma_1(t), \ldots, \gamma_n(t))$, the subunit condition means precisely that, for almost every $t$, we have

$$\left(\sum \dot\gamma_j(t)\xi_j\right)^2 \leq \sum a_{ij}(\gamma(t))\xi_i\xi_j$$

for any $\xi \in \mathbb{R}^n$, assuming the coefficient matrix $(a_{ij})$ to be symmetric.

The principal result of this paper is that these hypotheses do not suffice to guarantee global $C^\infty$ regularity.

**Theorem 1.1.** *There exists a selfadjoint second-order partial differential operator $L$ with $C^\infty$ coefficients and nonnegative principal symbol, defined on a smooth compact two-dimensional manifold without boundary, that is generically elliptic in the sense (1.1), has closed range on $L^2$ in the sense (1.2), and satisfies the reachability hypothesis (1.3), yet is not globally regular in $C^\infty$.*

Our example is closely related to one discovered recently [C] in the context of the $\bar\partial$-Neumann problem in several complex variables, where global regularity is a natural and important issue. The basic method of the present paper is that of [C], but some refinements are introduced here, in Propositions 2.3, 2.5 and 2.6. Our theorem may be interpreted in terms of the long time behavior of the diffusion generated by $L$. See Proposition 8.3.

## 2. Preliminaries and Results

Let $M$ be a compact, connected, two-dimensional $C^\infty$ manifold without boundary, and assume that, for some $\delta_0 > 0$, $M$ contains an open set $V_0$ diffeomorphic to

$$\{(x, t) \in \mathbb{R}^2 : |x| < 1 + \delta_0 \text{ and } |t| < \delta_0\}.$$

Throughout the discussion we consider $V_0$ to be identified with this subset of $\mathbb{R}^2$ via some fixed diffeomorphism. Assume that $M$ is equipped with a measure $\sigma$, having a nonvanishing $C^\infty$ density with respect to Lebesgue measure in any local coordinate system, and assume that $d\sigma \equiv dx\, dt$ in $V_0$. Let $L$ be a second-order partial differential operator on $M$, having real, $C^\infty$ coefficients and nonnegative principal symbol, that is selfadjoint in $L^2(M, \sigma)$. Assume that $L$ is elliptic at every point of $M \setminus J$, where $J = [-1, 1] \times \{0\} \subset V_0$, and that in $V_0$ we have

$$L = -\partial_x^2 - \partial_t \circ a(x,t)^2 \circ \partial_t + b(x,t),$$

where $a, b \in C^\infty$ are real-valued,

(2.1) $$a(x, 0) \equiv 0 \quad \text{for all } |x| \leq 1,$$

and

(2.2) $$\frac{\partial a}{\partial t}(x, 0) \neq 0 \quad \text{for all } |x| \leq 1$$

(that is, $\partial a/\partial t$ vanishes nowhere on $J$). The ellipticity hypothesis on $M \setminus J$ implies in particular that

(2.3) $$\{(x, t) \in V_0 : a(x, t) = 0\} = J.$$

We assume also that

(2.4) $$\|u\|_{L^2(M)} \leq C \|Lu\|_{L^2(M)} \quad \text{for all } u \in C^2(M),$$

and that $L$ is a bijection of its domain with $L^2(M)$. Thus $L^{-1}$ is a well-defined bounded linear operator on $L^2(M)$, whose range is the domain of $L$. In Section 4 we will show that such operators $L$ do exist; the only hypothesis that is not immediate is (2.4), and it can be ensured by choosing the coefficient $b$ to be sufficiently large.

Such an operator satisfies the reachability hypothesis (1.3). For outside $J$, any tangent vector is a scalar multiple of a subunit vector, so any two points in $M \setminus J$ may be joined via a subunit path, since $M \setminus J$ is connected. In $V_0$, $\partial_x$ is a subunit vector field, so any point of $J$ can be joined to a point in $V_0 \setminus J$ via a path $\gamma(s) = (c + s, 0)$.

Any such operator $L$ furnishes the desired example of global non-regularity.

**Theorem 2.1.** *Suppose that $L$ and $M$ are as described. Then $L$ is not globally regular in $C^\infty(M)$.*

A fundamental issue in this type of problem is that of exact regularity. Denote by $H^s(M)$ the Sobolev space of all functions in $L^2(M)$ possessing $s$ derivatives in $L^2(M)$, in the usual sense, for each nonnegative real $s$.

**Definition.** $L$ is said to be *exactly regular in $H^s$* if $u \in H^s$ for every $u \in L^2$ such that $Lu \in H^s(M)$, and if there is an a priori inequality $\|u\|_{H^s} \leq C_s \|Lu\|_{H^s}$ valid for all such $u$. We say that $L$ is *exactly regular* if it is exactly regular in $H^s$ for every $s \geq 0$.

Exact regularity implies global $C^\infty$ regularity, since $C^\infty$ is the intersection of all the Sobolev spaces. The examples $\partial_{x_1} + \alpha\partial_{x_2}$ of Section 1 are never exactly regular, yet are globally regular in $C^\infty$ for generic $\alpha$. No examples are known to this author, however, of second-order operators $L$ as described in Section 1, satisfying (1.1), (1.2) and (1.3), that are globally regular in $C^\infty$ yet are not exactly regular.

The proof of Theorem 2.1 rests on two results which are almost mutually contradictory. Assume that $L$ is as described above, satisfying conditions (2.1)–(2.4).

**Proposition 2.2.** *There exists $\hat{s} < \infty$ such that, for every $s \geq \hat{s}$, $L$ fails to be exactly regular in $\hat{s}$.*

**Proposition 2.3.** *There exists a discrete set $\Sigma \subset \mathbb{R}^+$ such that, for each $s \notin \Sigma$, there exists $C_s < \infty$ such that*

$$\|u\|_{H^s} \leq C_s \|Lu\|_{H^s} \tag{2.5}$$

*for every $u \in H^s(M)$ such that $Lu \in H^s(M)$.*

From these two propositions, Theorem 2.1 follows at once. Suppose that $L$ were globally regular in $C^\infty$. Fix any $s \geq \hat{s}$ that does not belong to $\Sigma$. Then, for any $f \in C^\infty(M)$, we have $L^{-1}f \in C^\infty(M)$, and hence $\|L^{-1}f\|_{H^s} \leq C_s\|f\|_{H^s}$ by (2.5). Since $L^{-1}$ is a bounded linear operator from $L^2(M)$ into $L^2(M)$ and $C^\infty$ is dense in $L^2$, it follows that $L^{-1}$ maps $H^s$ into $H^s$ and does so boundedly; in other words, $L$ is exactly regular in $H^s$. This contradicts Proposition 2.2. We have demonstrated the following:

**Corollary 2.4.** *Under the hypotheses of Theorem 2.1,*

$$\{f \in C^\infty(M) : L^{-1}f \in C^\infty(M)\}$$

*is not dense in $L^2(M)$.*

The next two results clarify the situation a bit. It will be shown that $\Sigma$ is discrete as a subset of $[0, \infty)$, not merely of $\mathbb{R}^+$. Write $\Sigma = \{0 < s_0 < s_1 < \cdots\}$.

**Proposition 2.5.** *If $r < s_0$, $u \in L^2$ and $Lu \in H^r$, then $u \in H^r$. If $s_j < s < r < s_{j+1}$, $u \in H^s$ and $Lu \in H^r$, then $u \in H^r$.*

**Proposition 2.6.** *For each $s \notin \Sigma$, the set $\{f \in H^s : L^{-1}f \in H^s\}$ is a closed subspace of $H^s$ having finite codimension in $H^s$.*

We believe this codimension to be equal to $2(j+1)$ if $s_j < s < s_{j+1}$.

## 3. A Model

In this section we analyze a boundary value problem that arises from the operator $L$ of Section 2 by a natural limiting process. Write

$$a(x,t) = \alpha(x)t + O(t^2),$$
$$b(x,t) = \beta(x) + O(t)$$

for $|x| \leq 1$. For each $s \in \mathbb{R}$ define

(3.1) $$\mathcal{L}_s = -\partial_x^2 - \alpha^2(x)(t\partial_t + s + 1)(t\partial_t + s) + \beta(x).$$

Consider the Dirichlet problem

(3.2) $$\begin{cases} \mathcal{L}_s u = f & \text{on } S, \\ u = 0 & \text{on } \partial S, \end{cases}$$

where $S = [-1,1] \times \mathbb{R}$ and $\partial S = \{-1,1\} \times \mathbb{R}$. We will often deal with solutions $u$ of (3.2) such that both $u$ and $\partial_x u$ belong to $L^2(S)$. This restriction on $\partial_x u$ suffices to ensure that $u(\pm 1, t)$ is well-defined as a function in $L^2(\mathbb{R})$, while the growth restriction on $u$ should be viewed as a boundary condition at infinity.

**Lemma 3.1.** *There exists a nonempty discrete set $\Sigma \subset [0, \infty)$, not containing $0$, having all of the following properties.*

  (i) *For each $s \notin \Sigma$, for each $f \in L^2(S)$ there exists a unique $u \in L^2(S)$ such that $\partial_x u \in L^2(S)$, $u(\pm 1, t) = 0$ for almost every $t \in \mathbb{R}$, and $\mathcal{L}_s u = f$ on $S$.*
  (ii) *If $s, f, u$ are as in (i) and if $f = f_1 + t\partial_t f_2 + (t\partial_t)^2 f_3$ where each $f_j \in L^2(S)$, then*
$$\|u\|_{L^2(S)} \leq C_s \sum \|f_j\|_{L^2(S)}.$$
  (iii) *For each $s \in \Sigma$ there exist $\gamma \in \mathbb{R}$ satisfying $\mathrm{Re}(\gamma) = s$ and $g \in C^\infty[-1,1]$ satisfying $g(\pm 1) = 0$, such that*
$$\mathcal{L}_0(g(x) t^{\gamma - \frac{1}{2}} \chi_{\mathbb{R}^+}(t)) \equiv 0.$$
  (iv) *A nonnegative real number $s$ belongs to $\Sigma$ if and only if there exists $u \in L^2(S)$, with compact support, such that $\mathcal{L}_0 u \in C^\infty$, $u(\pm 1, t) \equiv 0$, $u \in H^r$ for all $0 \leq r < s$, yet $u \notin H^s$.* □

Here $\chi_{\mathbb{R}^+}$ denotes the function that is identically $1$ for $t > 0$ and identically $0$ for $t < 0$. Condition (iv) follows immediately from (iii) by setting $u(x,t) = g(x) t^{\gamma - \frac{1}{2}} \chi_{\mathbb{R}^+}(t) \eta(t)$, where $\eta \in C_0^\infty(\mathbb{R})$ is identically $1$ in some neighborhood of $0$.

*Proof.* Define the partial Mellin transform by
$$\hat{f}(x,\tau) = \int_0^\infty f(x,t) t^{-i\tau} t^{-1}\, dt.$$

Set $S^\pm = [-1,1] \times \mathbb{R}^\pm$. We have $f \in L^2(S^+)$ with respect to Lebesgue measure if and only if $t^{1/2} f \in L^2(S^+, t^{-1}\, dt\, dx)$. In that event $\hat{f}(x, \tau + \tfrac{i}{2})$ is well-defined as a function in $L^2([-1,1] \times \mathbb{R})$, and

$$\|f\|_{L^2(S^+)}^2 = c_0 \iint_{[-1,1]\times\mathbb{R}} |\hat{f}(x, \tau + \tfrac{i}{2})|^2\, d\tau\, dx \tag{3.3}$$

for some finite constant $c_0$. If, moreover, $\mathcal{L}_s f \in L^2(S^+)$, then

$$(\mathcal{L}_s f)\hat{\,}(x, \tau + \tfrac{i}{2}) = H_z \hat{f}(x, \tau + \tfrac{i}{2}) \tag{3.4}$$

as functions in $L^2[-1,1]$, for almost every $\tau \in \mathbb{R}$, where
$$H_z = -\partial_x^2 - z(z+1)\alpha^2(x) + \beta(x)$$
and
$$z = s - \tfrac{1}{2} + i\tau,$$
as follows from the identity $(t\partial_t f)\hat{\,}(x,\gamma) = i\gamma \hat{f}(x,\gamma)$.

Set $\tilde{H}_w = -\partial_x^2 - w\alpha^2(x) + \beta(x)$ and define $\Sigma_0 \subset \mathbb{C}$ to be the set of all $w$ for which there exists $0 \neq g \in L^2[-1,1]$ satisfying $\tilde{H}_w g = 0$ and $g(\pm 1) = 0$. Both coefficients $\alpha, \beta$ are assumed real-valued, so that for any such $g$,
$$\operatorname{Im} \int_{-1}^1 \tilde{H}_w g \cdot \bar{g} = -\operatorname{Im}(w) \cdot \int_{-1}^1 \alpha^2 \cdot |g|^2.$$

Hypothesis (2.2) asserts that $\alpha$ vanishes nowhere on $[-1,1]$, so $\tilde{H}_w g = 0$ implies $\operatorname{Im}(w) = 0$. Thus $\Sigma_0 \subset \mathbb{R}$, and by examining instead the real part of $\int \tilde{H}_w g \cdot \bar{g}$ one finds in the same way that $\Sigma_0$ is bounded below.

$\Sigma_0$ is discrete, nonempty, and infinite. This is proved by considering the unique solution $g_w$ of $\tilde{H}_w g = 0$ satisfying $g_w(-1) = 0$ and $g_w'(-1) = 1$. The real-valued function that maps $w \in [0, \infty)$ to $g_w(1)$ must change sign infinitely many times, by a comparison argument, since $\alpha^2$ does not vanish.

Define
(3.5)
$$\Sigma = \left\{ s \in [0,\infty) : \text{there exists } t \in \mathbb{R} \text{ with } z(z+1) \in \Sigma_0, \text{ where } z = s - \tfrac{1}{2} + i\tau. \right\}$$

The set of all $z \in \mathbb{C}$ such that $z(z+1) \in \Sigma_0$ is discrete, and at most finitely many of its elements belong to $\mathbb{C} \setminus \mathbb{R}$. Therefore $\Sigma$ is also discrete. It will be shown later, in Corollary 5.2, that $0 \notin \Sigma$.

Conclusion (iii) follows directly from the definition of $\Sigma$, and has already been shown to imply (iv). Conclusions (i) and (ii) will be deduced from the following lemma, whose proof is nearly identical to that of Lemma 4 of [C] and is therefore omitted. Write $\langle \tau \rangle = (1 + |\tau|^2)^{1/2}$.

**Lemma 3.2.** *For each $s \in [0, \infty) \setminus \Sigma$ there exists $C_s < \infty$ such that, for every $\tau \in \mathbb{R}$, for every $g \in L^2[-1, 1]$ such that $H_{s-\frac{1}{2}+i\tau} g \in L^2[-1, 1]$, we have*

$$\|g\|_{L^2[-1,1]} + \langle \tau \rangle^{-1} \|\partial_x g\|_{L^2[-1,1]}$$
$$\leq C_s \langle \tau \rangle^{-2} \|H_{s-\frac{1}{2}+i\tau} g\|_{L^2[-1,1]} + C_s \langle \tau \rangle^{-1/2} (|g(1)| + |g(-1)|).$$

For fixed $\tau$ the conclusion follows directly from the condition that the nullspace of $H_{s-\frac{1}{2}+i\tau}$, with Dirichlet boundary conditions, is trivial. The bounds for $|\tau|$ large are obtained by examining $\langle H_{s-\frac{1}{2}+i\tau} g, g \rangle$ and integrating by parts. $\square$

To prove conclusion (i) of Lemma 3.1, let $s \notin \Sigma$ and $f \in L^2(S)$ be given. Then $\hat{f}(x, \tau + \frac{i}{2})$ is well-defined as an element of $L^2[-1, 1]$ for almost every $\tau \in \mathbb{R}$, and for all such $\tau$ we define $g(x, \tau)$ to be the unique solution of $H_{s-\frac{1}{2}+i\tau} g(x, \tau) = \hat{f}(x, \tau + \frac{i}{2})$ satisfying $g(\pm 1, \tau) = 0$. Define $u \in L^2(S^+)$ by $\hat{u}(x, \tau + \frac{i}{2}) = g(x, \tau)$ for almost every $\tau \in \mathbb{R}$. Since

$$\int_{-\infty}^{\infty} \int_{-1}^{1} |g(x, \tau)|^2 \, dx \, d\tau \leq C \int_{-\infty}^{\infty} \langle \tau \rangle^{-4} \int_{-1}^{1} |\hat{f}(x, \tau + \tfrac{i}{2})|^2 \, dx \, d\tau \leq C \|f\|_{L^2(S^+)}^2,$$

by Lemma 3.2, there does exist a unique $u \in L^2(S^+)$ with this prescribed partial Mellin transform. Since $(\partial_x u)\hat{}(x, \tau + \frac{i}{2}) = \partial_x g(x, \tau)$, the same reasoning shows that $\partial_x u \in L^2(S^+)$. Thus $\mathcal{L}_s u = f$ in $S^+$, by (3.4), and $u(\pm 1, t) = 0$ for almost all $t > 0$, by the choice of $g$.

A solution $u \in L^2(S^-)$ satisfying $\mathcal{L}_s u = f$ in $S^-$ may be found by the same procedure, after replacing $t$ by $-t$. Thus we obtain $u \in L^2(S)$ such that $\mathcal{L}_s u = f$, in the sense of distributions, on $(-1, 1) \times (0, \infty)$ and on $(-1, 1) \times (-\infty, 0)$. Moreover, $t\partial_t u \in L^2([-1, 1] \times \mathbb{R} \setminus \{0\})$, and likewise for $(t\partial_t)^2 u$, again by the bounds of Lemma 3.2. In order to show that $\mathcal{L}_s u = f$, in the sense of distributions, in $(-1, 1) \times \mathbb{R}$, it suffices to verify that $t\partial_t u \in L^2(S)$ and $(t\partial_t)^2 u \in L^2(S)$, in the sense of distributions. The variable $x$ plays no role then, and the desired result follows from the following remark.

**Lemma 3.3.** *If $h \in L^2(\mathbb{R})$ and $t\partial_t h \in L^2(\mathbb{R} \setminus \{0\})$ in the sense of distributions, then $t\partial_t h \in L^2(\mathbb{R})$.*

*Proof.* For any $\varphi \in C_0^\infty(\mathbb{R})$, for small $\varepsilon > 0$ we have
$$\int_0^\infty \big(\partial_t t\varphi(t)\big) h(t)\, dt = \int_\varepsilon^\infty \big(\partial_t t\varphi(t)\big) h(t)\, dt + O(\varepsilon^{1/2})$$
because $h \in L^2$. This is
$$\varepsilon\varphi(\varepsilon)h(\varepsilon) - \int_\varepsilon^\infty \varphi(t) \cdot t\partial_t h(t)\, dt + O(\varepsilon^{1/2}) = \varepsilon\varphi(\varepsilon)h(\varepsilon) + O(\varepsilon^{1/2}) - \int_0^\infty \varphi \cdot t\partial_t h,$$
since $t\partial_t h \in L^2$ and hence $h$ is continuous on $\mathbb{R}^+$. For any small $\delta > 0$, we have
$$\int_0^\delta \varepsilon^2 h^2(\varepsilon)\, d\varepsilon \leq \delta^2 \|h\|_{L^2}^2,$$
so there exists $\varepsilon \in (0, \delta]$ such that $|\varepsilon^2 h^2(\varepsilon)| \leq \delta \|h\|_{L^2}^2$. Hence there exists a sequence $\varepsilon_j \to 0^+$ such that $\varepsilon_j \varphi(\varepsilon_j) h(\varepsilon_j) \to 0$. Thus $\int_0^\infty \partial_t t\varphi \cdot h\, dt = -\int_0^\infty \varphi \cdot t\partial_t h\, dt$. The same reasoning applies on $(-\infty, 0]$. $\square$

Thus we have proved the existence statement in (i). To prove uniqueness, suppose that $u \in L^2(S)$, $\partial_x u \in L^2(S)$, $\mathcal{L}_s u = 0$, and $u(\pm 1, t) = 0$ for almost every $t$. Then $H_{s-\frac{1}{2}+i\tau} \hat{u}(x, \tau + \frac{i}{2}) \equiv 0$ and $\hat{u}(\pm 1, \tau + \frac{i}{2}) = 0$ for almost every $\tau$. But, by the definition of $\Sigma$, the operator $H_{s-\frac{1}{2}+i\tau}$ has trivial nullspace with Dirichlet boundary conditions, for every $\tau \in \mathbb{R}$. Thus $\hat{u}(x, \tau + \frac{i}{2}) \equiv 0$ for almost every $\tau \in \mathbb{R}$, so $u = 0$ in $L^2(S^+)$. Existence in $S^-$ is proved in the same way, using the substitution $t \mapsto -t$.

Conclusion (ii) follows in the same way by exploiting the factor of $\langle \tau \rangle^{-2}$ on the right-hand side in Lemma 3.2, since
$$((t\partial_t)^p f)\hat{\,}(x, \tau + \tfrac{i}{2}) = (i\tau - \tfrac{1}{2})^p \hat{f}(x, \tau + \tfrac{i}{2}).$$
If $\mathcal{L}_s u = \sum_{j \leq 2} (t\partial_t)^j f_j$, then
$$\|u\|_{L^2(S^+)}^2 = C \int_{-1}^1 \int_\mathbb{R} |\hat{u}(x, \tau+\tfrac{i}{2})|^2\, d\tau\, dx$$
$$\leq C \int_\mathbb{R} \langle \tau \rangle^{-1} |\hat{u}(\pm 1, \tau+\tfrac{i}{2})|^2\, d\tau + C \int_\mathbb{R} \int_{-1}^1 \langle \tau \rangle^{-4} \Big|\sum_j ((t\partial_t)^j f_j)\hat{\,}(x, \tau+\tfrac{i}{2})\Big|^2\, dx\, d\tau$$
$$\leq C\|u(\pm 1, \cdot)\|_{L^2(\mathbb{R})}^2 + C \int_\mathbb{R} \int_{-1}^1 \sum_j |\hat{f}_j(x, \tau+\tfrac{i}{2})|^2\, dx\, d\tau$$
$$\leq C\|u(\pm 1, \cdot)\|_{L^2(\mathbb{R})}^2 + C \sum_j \|f_j\|_{L^2(S)}^2.$$

The region $S^-$ is treated in the same way. In each line we are summing over both choices of the $\pm$ sign. $\square$

## 4. $L^2$ Theory

Fix any relatively compact neighborhood $V \subset V_0$ of $J$. For $u \in C^2(M)$, we have

$$(4.1) \qquad \int_{M \setminus V} |\nabla^2 u|^2 \leq C\|Lu\|^2 + C\|u\|^2 \leq C\|Lu\|^2$$

because $L$ is elliptic on $M \setminus J$, and $\|u\| \leq C\|Lu\|$. The constant $C$ in the first inequality does depend on $V$. Since $L = -\partial_x^2 - \partial_t a^2 \partial_t + b$ in $V_0$, integration by parts gives

$$\int_M Lu \cdot u \, d\sigma = \int_V (|\partial_x u|^2 + |a\partial_t u|^2 + b|u|^2) + O(\|\nabla u\|_{L^2(M\setminus V)}^2 + \|u\|_{L^2}^2).$$

Thus $\|\partial_x u\|_{L^2(V)} \leq C\|Lu\| + C\|u\|$, so we conclude that

$$(4.2) \qquad \|\partial_x u\|_{L^2(V)} \leq C\|Lu\| \quad \text{for all } u \in C^2(M),$$

where $\|\cdot\|$ with no subscript denotes the norm in $L^2(M, \sigma)$.

Consider $L$ as an unbounded, densely defined linear operator mapping $L^2(M)$ to $L^2(M)$. Define $\mathrm{Domain}(L)$ to be the closure of $C^2(M)$ under the graph norm $\|u\| + \|Lu\|$, which by hypothesis is equivalent to $\|Lu\|$. From (4.2) it then follows that, for all $u \in \mathrm{Domain}(L)$, we have $\partial_x u \in L^2(V)$, with norm bounded by $C\|Lu\|$.

**Lemma 4.1.** *If $u$ and $Lu$ belong to $L^2$ in the sense of distributions, then $u \in \mathrm{Domain}(L)$.*

*Proof.* Since $L$ is elliptic on $M \setminus J$, we have $u \in H^2_{\mathrm{loc}}$ on $M \setminus J$. Fix a cutoff function $\eta \in C^\infty$ that is supported in $V$ and is identically equal to 1 in some neighborhood of $J$. Since $u \in H^2$ in a neighborhood of the support of $1-\eta$, we have $(1-\eta)u \in \mathrm{Domain}(L)$, so it suffices to examine $\eta u = g$.

Identify $V$ with a subset of $\{(x,t) \in \mathbb{R}^2 : |x| < 1 + \delta_0, |t| < \delta_0\}$. Fix $\varphi \in C_0^\infty(\mathbb{R}^2)$ with $\int \varphi = 1$, and set $\varphi_\varepsilon(x,t) = \varepsilon^{-2} \varphi(\varepsilon^{-1}x, \varepsilon^{-1}t)$ and $P_\varepsilon f = f * \varphi_\varepsilon$ for any $f \in L^2(\mathbb{R}^2)$.

We claim that $\partial_x g \in L^2$ and $a\partial_t g \in L^2$ in the sense of distributions. This will follow if we prove that $\|\partial_x P_\varepsilon g\| + \|a\partial_t P_\varepsilon g\|$ is bounded above as $\varepsilon \to 0$. Write

$$\|\partial_x P_\varepsilon g\|^2 + \|a\partial_t P_\varepsilon g\|^2 + \int b(P_\varepsilon g)^2 = \langle g, P_\varepsilon^* L P_\varepsilon g \rangle$$
$$= \langle g, P_\varepsilon^* P_\varepsilon L g \rangle + \langle g, P_\varepsilon^*[L, P_\varepsilon]g \rangle.$$

Now $\partial_x$ commutes with $P_\varepsilon$, so

$$[L, P_\varepsilon] = -(\partial_t \circ a)[a\partial_t, P_\varepsilon] - [\partial_t \circ a, P_\varepsilon]a\partial_t + [b, P_\varepsilon]$$
$$= -(\partial_t \circ a)[a\partial_t, P_\varepsilon] - a\partial_t \circ [\partial_t \circ a, P_\varepsilon] - [[\partial_t \circ a, P_\varepsilon], a\partial_t] + [b, P_\varepsilon]$$
$$= -\partial_t \circ aB_\varepsilon^1 + B_\varepsilon^2,$$

where the $B_\varepsilon^j$ are pseudodifferential operators in the class $S_{1,0}^0$, uniformly as $\varepsilon \to 0$. Therefore

$$\|\partial_x P_\varepsilon g\|^2 + \|a\partial_t P_\varepsilon g\|^2 \leq C\|g\| \cdot \|Lg\| + C\|g\|^2 + |\langle a\partial_t P_\varepsilon g, B_\varepsilon^1 g\rangle|$$
$$\leq C\|g\|^2 + C\|Lg\|^2 + C\|g\| \cdot \|a\partial_t P_\varepsilon g\|.$$

Since all quantities in this inequality are known to be finite, the claim is proved.

To prove that $g \in \mathrm{Domain}(L)$, it suffices to show that $\|LP_\varepsilon g - Lg\| \to 0$ as $\varepsilon \to 0$, since $P_\varepsilon g \in C^\infty$ for all $\varepsilon > 0$. Since $P_\varepsilon Lg \to Lg$ in $L^2$ norm, it suffices to show that $\|[L, P_\varepsilon]g\| \to 0$. Now $[L, P_\varepsilon] = B_\varepsilon^1 \circ a\partial_t + B_\varepsilon^2$, where the $B_\varepsilon^j$ are again pseudodifferential operators in $S_{1,0}^0$, uniformly in $\varepsilon$. But these have the additional property that, on any compact subset of $T^*\mathbb{R}^2$, their principal symbols tend uniformly to 0, along with all of their derivatives, as $\varepsilon \to 0$. Therefore, since $a\partial_t g \in L^2$ and $g \in L^2$, it follows that $\|B_\varepsilon^1 a\partial_t g\| + \|B_\varepsilon^2 g\| \to 0$. □

**Corollary 4.2.** *$L$ is a bijection of* $\mathrm{Domain}(L)$ *onto* $L^2(M)$.

*Proof.* $L$ is injective and has closed range, by the hypothesis $\|u\| \leq C\|Lu\|$ for $u \in C^2$ and the definition of domain. If $f \in L^2$ is orthogonal to the range of $L$, then $L^*f = 0$ in the sense of distributions, where $L^*$ is the formal adjoint of $L$. But $L^* = L$, so $Lf = 0$, so $f \in \mathrm{Domain}(L)$ by Lemma 4.1, so $f = 0$. □

Observe now that there do exist operators $L$ satisfying hypothesis (2.1)–(2.4). There clearly exist operators $\tilde{L}$ satisfying (2.1)–(2.3). Integration by parts gives $\langle \tilde{L}u, u\rangle \geq -C\|u\|^2$ for all $u \in C^2$, because $\tilde{L}$ is symmetric and has nonnegative principal symbol. Then define $L = \tilde{L} + C + 1$.

## 5. Exact Regularity Does Not Hold

**Lemma 5.1.** *Suppose that $s \in [0, \infty]$ and $L$ is exactly regular on $H^s(M)$. Then, for each $f \in C_0^\infty((-1, 1) \times \mathbb{R})$, there exists $u \in L^2(S)$ such that $\partial_x u \in L^2(S)$, $|\partial_t|^s u \in L^2(S)$, and*

$$\|\partial_x u\|_{L^2(S)} + \||\partial_t|^s u\|_{L^2(S)} \leq C\|f\|_{H^s(S)},$$

*satisfying $u = 0$ almost everywhere on $\{\pm 1\} \times \mathbb{R}$ and $\mathcal{L}_0 u = f$.*

Here $|\partial_t|^s$ is defined in terms of the partial Fourier transform in the $t$ variable on $[-1,1] \times \mathbb{R}$.

*Proof.* Since $f \in C_0^\infty((-1,1) \times \mathbb{R})$, we may regard it as a function in $C_0^\infty(\mathbb{R}^2)$. For each small $\varepsilon > 0$, set
$$f_\varepsilon(x,t) = f(x, \varepsilon^{-1}t).$$
The function $f_\varepsilon$ is supported in $(-1,1) \times \{|t| \leq C\varepsilon\}$, so may be regarded as a function in $C^\infty(M)$, supported in $V$, for all sufficiently small $\varepsilon$. Then
$$\|f_\varepsilon\|_{H^s(M)} \leq C\varepsilon^{\frac{1}{2}-s}\|f\|_{H^s(S)}.$$
This follows by direct calculation for nonnegative integers $s$, then for general $s$ by interpolation. Set
$$v_\varepsilon = L^{-1}f_\varepsilon \quad \text{on } M$$
and
$$u_\varepsilon(x,t) = v_\varepsilon(x, \varepsilon t) \cdot \eta(\varepsilon t) \quad \text{on } S,$$
where $\eta \in C_0^\infty(\mathbb{R})$ is identically equal to 1 in some neighborhood of 0. We have
$$\|u_\varepsilon\|_{L^2(S)} \leq C\varepsilon^{-1/2}\|v_\varepsilon\|_{L^2(M)} \leq C\varepsilon^{-1/2}\|f_\varepsilon\|_{L^2(M)} \leq C\|f\|_{L^2(S)}.$$
Likewise,
$$\||\partial_t|^s u_\varepsilon\|_{L^2(S)} \leq C\|v_\varepsilon(x, \varepsilon t)\|_{L^2(S)} + c\||\partial_t|^s(v_\varepsilon(x, \varepsilon t))\|_{L^2(S)}$$
$$\leq C\|f\|_{L^2(S)} + C\varepsilon^{s-\frac{1}{2}}\|v_\varepsilon\|_{H^s(M)}$$
$$\leq C\|f\|_{L^2(S)} + C\varepsilon^{s-\frac{1}{2}}\|f_\varepsilon\|_{H^s(M)}$$
by the hypothesis of exact regularity on $H^s(M)$. Thus
$$\||\partial_t|^s u_\varepsilon\|_{L^2(S)} \leq C\|f\|_{H^s(S)}.$$
Since (4.2) holds for all functions in Domain($L$), we likewise have
$$\|\partial_x u_\varepsilon\|_{L^2(S)} \leq C\varepsilon^{-1/2}\|\partial_x v_\varepsilon\|_{L^2(M)} \leq C\varepsilon^{-1/2}\|f_\varepsilon\|_{L^2(M)} \leq C\|f\|_{L^2(S)}.$$

There exists a sequence $\varepsilon_j \to 0$ such that $u_{\varepsilon_j}$, $\partial_x u_{\varepsilon_j}$, and $|\partial_t|^s u_{\varepsilon_j}$ converge in the weak $*$ topology in $L^2(S)$. Defining $u \in L^2(S)$ to be the limit of one such subsequence, we have
$$\|u\|_{L^2(S)} + \|\partial_x u\|_{L^2(S)} \leq C\|f\|_{L^2(S)}$$
and
$$\||\partial_t|^s u\|_{L^2(S)} \leq C\|f\|_{H^s(S)}.$$
We claim that $\mathcal{L}_0 u = f$, in the sense of distributions on the interior of $S$. To verify this, consider any test function $\varphi \in C_0^\infty((-1,1) \times \mathbb{R})$. Then $\eta(\varepsilon t) \equiv 1$ on the support of $\varphi$ for all sufficiently small $\varepsilon$, so on the support of $\varphi$, $u$ is the weak limit of $v_{\varepsilon_j}(x, \varepsilon_j t)$. Now writing $\varphi^\varepsilon(x, r) = \varphi(x, \varepsilon^{-1} r)$,

$$\begin{aligned}
\langle u_\varepsilon, \mathcal{L}_0 \varphi \rangle &= \iint_S v_\varepsilon(x, \varepsilon t)(\mathcal{L}_0 \varphi)(x, t)\, dx\, dt \\
&= \varepsilon^{-1} \iint v_\varepsilon(x, r)(\mathcal{L}_0 \varphi)(x, \varepsilon^{-1} r)\, dx\, dr \\
&= \varepsilon^{-1} \iint f_\varepsilon(x, r) \varphi(x, \varepsilon^{-1} r)\, dx\, dr \\
&\quad + \varepsilon^{-1} \iint (\mathcal{L}_0 - \mathcal{L}) v_\varepsilon(x, r) \cdot \varphi^\varepsilon(x, r)\, dx\, dr \\
&= \iint f \varphi\, dx\, dt + \varepsilon^{-1} \iint v_\varepsilon(x, r)((\mathcal{L}_0 - \mathcal{L})\varphi^\varepsilon)(x, r)\, dx\, dr.
\end{aligned}$$

The difference $\mathcal{L}_0 - \mathcal{L}$ may be written in the $(x, r)$ coordinates as a sum of three terms $O(r)$, $O(r^2)\partial_r$, and $O(r^3)\partial_r^2$, where $O(r^n)$ denotes a $C^\infty$ function on $[-1, 1] \times \mathbb{R}$ that is $O(r^n)$ as $r \to 0$. Since $\varphi$ has compact support,

$$\|O(r) \cdot \varphi^\varepsilon\|_{L^2(S)} = O(\varepsilon)\|\varphi^\varepsilon\|_{L^2(S)} = O(\varepsilon^{3/2}).$$

Likewise,

$$\|O(r^2)\partial_r \varphi^\varepsilon\|_{L^2(S)} + \|O(r^3)\partial_r^2 \varphi^\varepsilon\|_{L^2(S)} = O(\varepsilon^{3/2}),$$

for fixed $\varphi \in C_0^\infty$. Thus as $\varepsilon \to 0$,

$$\begin{aligned}
\left| \varepsilon^{-1} \iint v_\varepsilon(x, r)((\mathcal{L}_0 - \mathcal{L})\varphi^\varepsilon)(x, r)\, dx\, dr \right| &\leq C\varepsilon^{-1} \|v_\varepsilon\|_{L^2(S)} \varepsilon^{3/2} \\
&\leq C\varepsilon^{-1} \varepsilon^{1/2} \|f\|_{L^2(S)} \cdot \varepsilon^{3/2} \\
&= O(\varepsilon).
\end{aligned}$$

Since $u_\varepsilon \to u$ in the weak * topology, we conclude that

$$\langle u, \mathcal{L}_0 \varphi \rangle = \langle f, \varphi \rangle$$

in the inner product for $L^2(S)$, so $\mathcal{L}_0 u = f$ in the sense of distributions.

The only conclusion of this lemma that may be unexpected is that the weak * limit $u$ automatically satisfies the Dirichlet boundary condition $u(\pm 1, t) = 0$ almost everywhere. This follows from the hypoellipticity of $L$ where $|x| > 1$, coupled with a strong form of the reachability hypothesis. To see this, fix a small parameter $\delta > 0$. Since $L$ is elliptic on $M \setminus J$ and $Lv_\varepsilon \equiv 0$ on $M \setminus \{(x, t) \in V : |x| \leq 1\}$, elliptic regularity gives a bound

$$|v_\varepsilon(1 + \delta, t)| \leq C_\delta \|v_\varepsilon\|_{L^2(M)} \leq C_\delta \varepsilon^{1/2}$$

for all $|t| \leq \delta_0$, uniformly in $\varepsilon$, where $C_\delta$ depends also on $f$. Thus, for any constant $A < \infty$, for all sufficiently small $\varepsilon$ we have

$$\left(\int_{|t|\leq \varepsilon A} v_\varepsilon(1,t)^2\, dt\right)^{1/2} \leq \left(\int_{|t|\leq \varepsilon A} |v_\varepsilon(1,t) - v_\varepsilon(1+\delta, t)|^2\, dt\right)^{1/2} + C_{\delta,A}\varepsilon$$

$$\leq \delta^{1/2}\left(\int_{|t|\leq \varepsilon A}\int_1^{1+\delta} |\partial_y v_\varepsilon(y,t)|^2\, dy\, dt\right)^{1/2} + C_{\delta,A}\varepsilon$$

$$\leq \delta^{1/2}\|\partial_x v_\varepsilon\|_{L^2(M)} + C_{\delta,A}\varepsilon$$

$$\leq C\delta^{1/2}\varepsilon^{1/2} + C_{\delta,A}\varepsilon.$$

Hence

$$\left(\int_{|r|\leq A} |u_\varepsilon(1,r)|^2\, dr\right)^{1/2} = \left(\int_{|r|\leq A} |v_\varepsilon(1,\varepsilon r)|^2\, dr\right)^{1/2} \leq C\delta^{1/2} + C_{\delta,A}\varepsilon^{1/2}.$$

Taking the limit as $\varepsilon \to 0$ gives

$$\int_{|r|\leq A} |u(1,r)|^2\, dr \leq C\delta,$$

for $\int_{|r|\leq A}|u_\varepsilon(1,r)|^2\, dr$ tends to the corresponding integral of $u$ as $\varepsilon \to 0$, because $u_\varepsilon \to u$ in the weak * sense on $S$ and $\partial_x u_\varepsilon$ is uniformly bounded in $L^2(S)$. Since $A$ and $\delta$ are arbitrary positive constants, $u(1,r) = 0$ for almost every $r \in \mathbb{R}$. The same reasoning applies to $u(-1,r)$. $\square$

**Corollary 5.2.** $0 \notin \Sigma$.

*Proof.* Applying Lemma 5.1 with $s = 0$, we find that for every $f \in L^2(S)$ there exists $u \in L^2(S)$ satisfying the Dirichlet boundary conditions and $\mathcal{L}_0 u = f$. Thus $H_{-\frac{1}{2}+i\tau}\hat{u}(x, \tau+\frac{i}{2}) = \hat{f}(x, \tau+\frac{i}{2})$ for almost every $\tau \in \mathbb{R}$, and $\hat{u}(\pm 1, \tau+\frac{i}{2}) = 0$. But for all but finitely many $\tau \in \mathbb{R}$, regardless of whether $0 \in \Sigma$ or not, $\hat{u}(\cdot, \tau+\frac{i}{2})$ is uniquely determined by this ordinary differential equation and boundary condition, since $\Sigma_0$ is discrete. Thus $u$ is uniquely determined by $f$, via the partial Mellin transform followed by $H^{-1}_{-\frac{1}{2}+i\tau}$, followed by the inverse Mellin transform. If $0 \in \Sigma$, there exists $\tau_0 \in \mathbb{R}$ such that the operator norm of $H^{-1}_{-\frac{1}{2}+i\tau}$ tends to infinity as $\tau \to \tau_0$. Hence, by the Plancherel identity for the Mellin transform, the operator $f \to u$ is not bounded in $L^2(S)$. Thus the assumption $0 \in \Sigma$ leads to a contradiction. $\square$

**Corollary 5.3.** *If $s > 0$ and $L$ is exactly regular in $H^s(M)$, then $s \notin \Sigma$.*

Since exact regularity in $H^s$ implies exact regularity in $H^r$ for all $0 \le r \le s$, by interpolation, the full conclusion is that, if
$$s_0 = \min_{s \in \Sigma} s,$$
then
> For all $s \ge s_0$, the operator $L$ is not exactly regular in $H^s(M)$.

Recall that $\Sigma$ has been shown to be nonempty.

*Proof of Corollary 5.3.* This follows from Lemma 5.1 in essentially the same way as does Corollary 5.2. Suppose that $L$ is exactly regular in $H^s(M)$, yet also that $s \in \Sigma$. By Lemma 3.1 there exists $f \in C_0^\infty(S)$ for which there exists $u \in L^2(S)$ satisfying $\mathcal{L}_0 u = f$ and $u(\pm 1, \cdot) \equiv 0$, but $u \notin H^s(S)$ and more precisely, by conclusion (iii) of that lemma, $|\partial_t|^s u \notin L^2(S)$. Lemma 5.1, on the other hand, yields a solution $\tilde{u} \in L^2(S)$ of the same boundary value problem, with $|\partial_t|^s \tilde{u} \in L^2(S)$ and $\partial_x \tilde{u} \in L^2(S)$. Since $0 \notin \Sigma$, the solution of the boundary value problem is unique in $L^2(S)$. Thus $u = \tilde{u}$, a contradiction. $\square$

The proof of Proposition 2.2 is now complete.

## 6. The Closed Range Estimate

In this section we establish the inequality

(6.1) $\qquad \|u\|_{H^s(M)} \le C_s \|Lu\|_{H^s(M)} \quad$ for all $u \in H^s$ such that $Lu \in H^s$,

for all $0 < s \notin \Sigma$. From this it follows that $L$ has closed range, as an unbounded densely defined operator from $H^s$ to $H^s$.

To begin, several preliminary reductions and estimates will be made, in which the hypothesis $s \notin \Sigma$ plays no role. Fix $s > 0$. Let $\gamma > 0$ be a small constant, to be chosen at the end of the section, and define
$$V_\gamma = \{(x,t) : |x| < 1+\gamma \text{ and } |t| < \gamma\}.$$
For any $f \in H^s(M)$, we have $L^{-1} f \in H^{s+2}$ in $M \setminus V_{\gamma/2}$, with norm bounded by $C(\gamma)\|f\|_{H^s}$, by elliptic regularity and the $L^2$ boundedness of $L^{-1}$. Fix $\eta \in C_0^\infty(V_\gamma)$ satisfying $\eta \equiv 1$ on $V_{\gamma/2}$. If $L^{-1} f \in H^s$ then $L(\eta L^{-1} f) = \eta f + [L, \eta] L^{-1} f$ belongs to $H^s$ with norm bounded by $C(\gamma)\|L^{-1} f\|_{H^s} + C(\gamma)\|f\|_{H^s}$, while
$$\|(1-\eta) L^{-1} f\|_{H^{s+2}} \le C(s,\gamma)\|f\|_{H^s}.$$
Thus it suffices to show that, for some small $\gamma > 0$, (6.1) holds for all $u$ supported in $V_\gamma$. Identifying $V_\gamma$ with a subset of $\mathbb{R}^2$, we may work henceforth in $\mathbb{R}^2$.

Let $(x, t, \xi, \tau)$ be coordinates in $T^*\mathbb{R}^2$, with $\xi$ dual to $x$ and $\tau$ to $t$. Define
$$\Gamma = \{(x, t, \xi, \tau) : |x| \leq 1,\ t = 0 \text{ and } \xi = 0\}.$$
$\Gamma$ is the set of all points where the principal symbol of $L$ vanishes. Therefore, for any classical pseudodifferential operator $Q \in S^0_{1,0}$ whose full symbol vanishes identically in some conic neighborhood of $\Gamma$ outside a bounded set, we have

(6.2) $$\|Qu\|_{H^{s+2}} \leq C_Q \|Lu\|_{H^s} + C_Q \|u\|_{L^2}$$

for all $u \in H^s$, supported in $V_{\delta_0}$, such that $Lu \in H^s$. The last term is of course redundant, being majorized by $C\|Lu\|_{L^2}$.

**Lemma 6.1.** *If $s \geq 0$, $u \in H^s$ is supported in $V_{\delta_0/2}$ and $Lu \in H^s$, then $\partial_x u$ and $a\partial_t u$ belong to $H^s$, and*
$$\|\partial_x u\|_{H^s} + \|a\partial_t u\|_{H^s} \leq C\|Lu\|_{H^s} + C\|u\|_{H^s}.$$

*Proof.* Fix $\varphi \in C_0^\infty(\mathbb{R}^2)$ satisfying $\varphi \equiv 1$ in a neighborhood of 0. Define the Fourier multiplier operators
$$(P_\varepsilon g)\hat{}\,(\xi, \tau) = \varphi(\varepsilon\xi, \varepsilon\tau)(1 + \xi^2 + \tau^2)^{s/2}\hat{g}(\xi, \tau),$$
mapping functions $g$ defined on $\mathbb{R}^2$ to functions defined on $\mathbb{R}^2$. Fix $\eta \in C_0^\infty(V_{\delta_0})$ that is identically equal to 1 on a neighborhood of $V_{\delta_0/2}$, and consider, for small $\varepsilon > 0$,
$$\|\partial_x \eta P_\varepsilon u\|^2 + \|a\partial_t \eta P_\varepsilon u\|^2 + \int b \cdot (\eta P_\varepsilon u)^2 = \langle \eta P_\varepsilon u, L\eta P_\varepsilon u\rangle$$
$$= \langle \eta P_\varepsilon u, \eta P_\varepsilon Lu\rangle + \langle \eta P_\varepsilon u, [L, \eta P_\varepsilon]u\rangle.$$

Now $[L, \eta P_\varepsilon] = \partial_x B_1 + a\partial_t B_2 + B_3$, where $B_1, B_2, B_3 \in S^s_{1,0}$ uniformly as $\varepsilon \to 0$. Thus
$$\|\partial_x \eta P_\varepsilon u\|^2 + \|a\partial_t \eta P_\varepsilon u\|^2 \leq C\|u\|_{H^s}\|Lu\|_{H^s} + C(\|\partial_x \eta P_\varepsilon u\| + \|a\partial_t \eta P_\varepsilon u\| + \|\eta P_\varepsilon u\|)\cdot\|u\|_{H^s}.$$
Therefore
$$\|\partial_x \eta P_\varepsilon u\| + \|a\partial_t \eta P_\varepsilon u\| \leq C\|Lu\|_{H^s} + C\|u\|_{H^s},$$
uniformly as $\varepsilon \to 0$. Invoking again the fact that the commutators of $\partial_x$ and $a\partial_t$ with $\eta P_\varepsilon$ belong to $S^s_{1,0}$ uniformly in $\varepsilon$, the lemma follows. $\square$

Define
$$(\Lambda^s g)\hat{}\,(\xi, \tau) = (1 + \xi^2 + \tau^2)^{s/2}\hat{g}(\xi, \tau)$$
for all $g$ in the Schwartz class, with $\hat{g}$ denoting here the Fourier transform. pseudodifferential operators in the class $S^0_{1,0}$, which are permitted to change from one occurrence of either symbol to the next, even within the same line. By a classical operator we mean one whose full symbol, in the Kohn–Nirenberg calculus, admits an asympotic expansion in terms homogeneous of integral degrees. Moreover $A$ denotes always an operator whose principal symbol $\sigma_0(A)$ vanishes identically on $\Gamma$.

The desired inequality (6.1), for all $u$ supported in $V_{\gamma/2}$, would follow via the substitution $g = \Lambda^s u$ from

$$\|g\| \leq C\|\Lambda^s L \Lambda^{-s} g\| + C\|g\|_{H^{-1}} \tag{6.3}$$

for all $g \in L^2$ supported in $V_\gamma$ such that $\Lambda^s L \Lambda^{-s} g \in L^2$, except that $L \circ \Lambda^{-s}$ is not defined since $\Lambda^{-s} g$ will in general not have compact support, nor will $g = \Lambda^s u$ have compact support. But fixing $h \in C_0^\infty(V_\gamma)$ identically equal to 1 on a neighborhood of the closure of $V_{\gamma/2}$, inequality (6.1) for all $u$ supported in $V_{\gamma/4}$ does follow from

$$\|g\| \leq C\|h\Lambda^s L h \Lambda^{-s} g\| + C\|g\|_{H^{-1}}, \tag{6.4}$$

by the substitution $g = h\Lambda^s u$, because $\Lambda^{\pm s}$ is pseudolocal. Define $L_s g = h\Lambda^s L h \Lambda^{-s} g$ for all $g$ supported in $V_{\gamma/2}$. Thus essentially

$$L_s = \Lambda^s L \Lambda^{-s}.$$

**Lemma 6.2.** *As an operator from functions supported on $V_{\gamma/2}$ to functions supported on $V_\gamma$, modulo an operator smoothing of infinite order, $L_s$ takes the form*

$$L_s = -\partial_x^2 - (\tilde{Y}_s + A)(Y_s + A) + \beta(x) + A, \tag{6.5}$$

*where $Y_s$ and $\tilde{Y}_s$ are real $C^\infty$ vector fields which, for $|x| \leq 1$, take the forms*

$$Y_s = \alpha(x)(t\partial_t + s) + O(t^2)\partial_t,$$
$$\tilde{Y}_s = \alpha(x)(t\partial_t + s + 1) + O(t^2)\partial_t.$$

*Moreover, $-\tilde{Y}_s$ is the formal adjoint of $Y_s$, modulo an operator in $S_{1,0}^0$.*

The terms $O(t^2)\partial_t$ depend on $s$, as do the operators $A$.

The proof is a straightforward calculation that is essentially identical to that in the second part of Section 3 of [C], so is omitted here. The infinitely smoothing operator by which $L_s$ differs from the operator in (6.5) may be neglected in establishing the inequality

$$\|u\| \leq C\|L_s u\| + C\|u\|_{H^{-1}} \quad \text{for all } u \text{ supported in } V_{\gamma/2}, \tag{6.6}$$

since its contribution is bounded by $C\|u\|_{H^{-1}}$ because it is smoothing.

From Lemma 6.1 it follows that $\partial_x u \in L^2$ and $a\partial_t u \in L^2$ if both $u$ and $L_s u$ belong to $L^2$, and

$$\|\partial_x u\| + \|a\partial_t u\| \leq C\|u\| + C\|L_s u\|. \tag{6.7}$$

**Lemma 6.3.** *For any small $\gamma > 0$ and any $u \in L^2$ supported in $V_\gamma$ such that $\partial_x u \in L^2$, we have*

$$\int |u(\pm 1, t)|^2 \, dt \leq C\gamma \|\partial_x u\|^2 \tag{6.8}$$

*and*

$$\iint_{1 < |x| < 1+\gamma} |u(x,t)|^2 \, dx \, dt \leq C\gamma^2 \|\partial_x u\|^2.$$

*Proof.* We have

$$|u(1,t)| = |u(1,t) - u(1+\gamma, t)| \leq \gamma^{1/2} \left( \int_1^{1+\gamma} |\partial_x u(r,t)|^2 \, dr \right)^{1/2}$$

for almost every $t \in \mathbb{R}$, and (6.8) follows directly. The proof of the second inequality is similar. $\square$

If $\gamma$ is chosen to be sufficiently small, depending on $s$, then combining the second inequality of the lemma with (6.7) gives

$$\|\partial_x u\| + \|u\| + \|a\partial_t u\| \leq C\|L_s u\| + C\|u\|_{H^{-1}} + C\|u\|_{L^2(S)}, \tag{6.9}$$

assuming that $u \in L^2$ is supported in $V_\gamma$ and that $L_s u \in L^2$, where $S = [-1, 1] \times \mathbb{R}$.

Recall the operators $\mathcal{L}_s = -\partial_x^2 - \alpha^2(x)(t\partial_t + s + 1)(t\partial_t + s) + \beta(x)$ introduced in Section 3. For $(x, t) \in S$ we may express

$$(L_s - \mathcal{L}_s)u = t\partial_t A_1 t\partial_t u + t\partial_t A_2 u + A_3 u$$

for certain $A_j \in S^0_{1,0}$ satisfying $\sigma_0(A_j)(x, 0, 0, \tau) \equiv 0$ for $|x| \leq 1$. This holds because multiplication by $t$ is an operator of the type $A$, as is $[t, A]\partial_t$ because $[t, A] \in S^{-1}_{1,0}$ has principal symbol $\sigma_{-1} = c \cdot \partial_\tau \sigma_0(A)$, which vanishes on $\Gamma$ because $\sigma_0(A)$ itself vanishes, and $\partial_\tau$ is tangent to $\Gamma$.

Assume now that $s \notin \Sigma$. Then there exists $w \in L^2(S^+)$ such that $\partial_x w \in L^2(S^+)$, $w(\pm 1, t) = 0$ for almost all $t$, $\mathcal{L}_s w = -(L_s - \mathcal{L}_s)u$, and

$$\|w\|_{L^2(S^+)} \leq C\|A_1 t\partial_t u\| + C\|A_2 u\| + C\|A_3 u\|. \tag{6.10}$$

Indeed, define $\hat{w}(x, \tau + \tfrac{i}{2})$, for almost every $\tau \in \mathbb{R}$, to be the unique solution of

$$H_{s-\frac{1}{2}+i\tau} \hat{w}(x, \tau + \tfrac{i}{2})$$
$$= i(\tau + \tfrac{i}{2})(A_1 t\partial_t u)\hat{\,}(x, \tau + \tfrac{i}{2}) + i(\tau + \tfrac{i}{2})(A_2 u)\hat{\,}(x, \tau + \tfrac{i}{2}) + (A_3 u)\hat{\,}(x, \tau + \tfrac{i}{2})$$

satisfying $\hat{w}(\pm 1, \tau + \tfrac{i}{2}) = 0$, and define $w$ to be the inverse Mellin transform of $\hat{w}$. Lemma 3.2 implies that $w$ and $\partial_x w$ belong to $L^2(S^+)$, and gives the inequality (6.10).

We thus have $\mathcal{L}_s(u-w) = L_s u + (\mathcal{L}_s - L_s)u - \mathcal{L}_s w = L_s u$ in $S^+$, $u - w = u$ on $\partial S^+$, and $\partial_x(u-w) \in L^2(S^+)$. Therefore, since $s \notin \Sigma$, Lemma 3.1 gives

$$\|u-w\|_{L^2(S^+)} \le C\|L_s u\|_{L^2(S^+)} + C\left(\int |u(\pm 1, t)|^2\, dt\right)^{1/2},$$

where we sum over both choices of the $\pm$ sign, and hence

$$\|u\|_{L^2(S^+)} \le C\|L_s u\| + C\left(\int |u(\pm 1, t)|^2\, dt\right)^{1/2} + C\|At\partial_t u\|_{L^2(S^+)} + C\|Au\|_{L^2(S^+)}.$$

The same reasoning may be repeated on $S^- = [-1,1] \times \mathbb{R}^-$ to obtain

$$\|u\|_{L^2(S)} \le C\|L_s u\| + C\left(\int |u(\pm 1, t)|^2\, dt\right)^{1/2} + C\|At\partial_t u\|_{L^2(S)} + C\|Au\|_{L^2(S)}.$$

Given any $A$ and any $\varepsilon > 0$, there exists $Q \in S^0_{1,0}$ whose full symbol vanishes identically on a conic neighborhood of $\Gamma$, such that $A = Q + R$ with $|\sigma_0(R)(x,t,\xi,\tau)| \le \varepsilon$ for all $(x,t,\xi,\tau) \in T^*\mathbb{R}^2$. Therefore, by a basic theorem on pseudodifferential operators,

$$\|Rg\|_{L^2} \le C\varepsilon\|g\|_{L^2} + C_\varepsilon\|g\|_{H^{-1}} \quad \text{for all } g \in L^2,$$

where $C$ is a universal constant. On the other hand, we have already seen that

$$\|Qu\|_{H^2} \le C_\varepsilon\|L_s u\|_{L^2} + C_\varepsilon\|u\|_{H^{-1}},$$

since $L_s$ is elliptic except on $\Gamma$. Thus

$$\|At\partial_t u\| + \|Au\| \le C_\varepsilon\|L_s u\| + C_\varepsilon\|u\|_{H^{-1}} + \varepsilon\|t\partial_t u\|_{L^2(S)} + \varepsilon\|u\|,$$

so

$$\|u\|_{L^2(S)} \le C\|L_s u\| + C\|u\|_{H^{-1}} + C\gamma^{1/2}\|\partial_x u\| + \varepsilon\|u\| + \varepsilon\|a\partial_t u\|,$$

where $C$ depends on $\varepsilon$. This, together with (6.9), gives

$$\|u\| + \|\partial_x u\| + \|a\partial_t u\| \le C\|L_s u\| + C\|u\|_{H^{-1}} + C\gamma^{1/2}\|\partial_x u\| + \varepsilon\|u\| + \varepsilon\|a\partial_t u\|_{L^2(S)}.$$

Since $u$, $\partial_x u$, and $a\partial_t u$ are already known to belong to $L^2$, and since $|a(x,t)| \ge c|t|$ for all $(x,t) \in [-1,1] \times (-\delta, \delta)$ for some $c, \delta > 0$, we may choose $\varepsilon, \gamma > 0$ so that the last three terms on the right-hand side may be absorbed into the left, yielding the desired inequality

$$\|u\| \le C\|L_s u\| + C\|u\|_{H^{-1}}$$

for all $u \in L^2$ supported in $V_\gamma$ satisfying $L_s u \in L^2$.

## 7. Limited Regularity

In this section we prove Propositions 2.5 and 2.6, which refine the preceding analysis. The main technique is a regularization process more delicate than that required in conventional hypoellipticity proofs.

*Proof of Proposition 2.5.* Supposing that $u \in H^s$ and $Lu \in H^r$, where either $0 \leq s < r < s_0$ or $s_j < s < r < s_{j+1}$ for some $j \geq 0$, we claim that $u \in H^r$. It suffices to prove this merely for some $r > s$, provided that $r - s$ is bounded below by a strictly positive quantity for all $s$ in any compact subset of $(s_j, s_{j+1})$ (or of $[0, s_0)$), since iterating such a result finitely many times yields the general case.

Let $s \notin \Sigma$ be given, and let $\delta > 0$ be a small number depending on $s$, to be determined. Fix a strictly positive function $h \in C^\infty(\mathbb{R}^2)$ satisfying $h(y) \equiv 1$ for $|y| \leq 1$, and $h(y) \equiv |y|^{-1}$ for all $|y| \geq 2$. For each small $\varepsilon > 0$ define the regularizing operator
$$(P_\varepsilon g)\hat{\ }(\xi, \tau) = h^\delta(\varepsilon\xi, \varepsilon\tau)\hat{g}(\xi, \tau).$$
Set $r = s + \delta$. Then $u \in H^s$ implies $P_\varepsilon u \in H^r$ for all $\varepsilon > 0$. Moreover
$$P_\varepsilon Lu = (P_\varepsilon L P_\varepsilon^{-1}) P_\varepsilon u.$$
Setting $g_\varepsilon = P_\varepsilon u$ and $L_{(\varepsilon)} = P_\varepsilon L P_\varepsilon^{-1}$, we have $g_\varepsilon \in H^r$ for all $\varepsilon > 0$, while $L_{(\varepsilon)} g_\varepsilon$ belongs to $H^r$ with a norm that is bounded uniformly as $\varepsilon \to 0$. If the analysis of Section 6 can be shown to apply to $L_{(\varepsilon)}$, uniformly as $\varepsilon \to 0$, we could conclude that $\|g_\varepsilon\|_{H^r}$ is bounded uniformly as $\varepsilon \to 0$, and consequently that $u \in H^r$.

As in Section 6, we may reduce to the case of functions supported in $V_\gamma$ for small $\gamma$. Now $P_\varepsilon \circ a\partial_t \circ P_\varepsilon^{-1} - a\partial_t = P_\varepsilon[a\partial_t, P_\varepsilon^{-1}]$ is an operator in $S^0_{1,0}$ that depends on both parameters $\varepsilon$ and $\delta$. Its principal symbol, for instance, is majorized by
$$Ch^\delta(\varepsilon\xi, \varepsilon\tau) \cdot |\tau| \cdot |\nabla_{\xi,\tau}[h^{-\delta}(\varepsilon\xi, \varepsilon\tau)]| \leq C\delta|\varepsilon\tau| \, |\nabla h/h|(\varepsilon\xi, \varepsilon\tau) \leq C\delta,$$
and in fact any $S^0_{1,0}$ seminorm of its full symbol is $O(\delta)$, uniformly in $\varepsilon$. Thus (in $V_\gamma$) we have
$$L - L_{(\varepsilon)} = B_1 \circ a\partial_t + B_2,$$
where $B_1, B_2 \in S^0_{1,0}$ are $O(\delta)$ in all $S^0_{1,0}$ seminorms, uniformly as $\varepsilon \to 0$. Therefore if $\delta$ is chosen sufficiently small, depending on the distance from $s$ to $\Sigma$, Proposition 2.3 yields for any $g$ supported in $V_\gamma$ a bound
$$\begin{aligned}\|g\|_{H^r} + \|a\partial_t g\|_{H^r} &\leq C\|Lg\|_{H^r} \\ &\leq C\|L_{(\varepsilon)} g\|_{H^r} + C\|B_1 a\partial_t g\|_{H^r} + C\|B_2 g\|_{H^r} \\ &\leq C\|L_{(\varepsilon)} g\|_{H^r} + C\delta\|a\partial_t g\|_{H^r} + C\delta\|g\|_{H^r}.\end{aligned}$$

If we can show that $a\partial_t g_\varepsilon$ belongs to $H^r$ then for sufficiently small $\delta$ we may absorb terms from the right hand side in this inequality to conclude that

$$\|g_\varepsilon\|_{H^r} \le C\|L_{(\varepsilon)}g_\varepsilon\|_{H^r} \le C < \infty,$$

uniformly as $\varepsilon \to 0$. Now

$$a\partial_t P_\varepsilon u = P_\varepsilon a\partial_t u + [a\partial_t, P_\varepsilon]u.$$

Since $u, Lu \in H^s$, we have $a\partial_t u \in H^s$ by Lemma 6.1, and therefore $P_\varepsilon a\partial_t u \in H^{s+\delta} = H^r$. The commutator $[a\partial_t, P_\varepsilon]$ is an operator of order $-\delta$, so maps $H^s$ to $H^r$, so $[a\partial_t, P_\varepsilon]u \in H^r$ as well. $\square$

*Proof of Proposition 2.6.* This proposition asserts that $\{f \in H^s : L^{-1}f \in H^s\} = W^s$ is a closed subspace of $H^s$ of finite codimension, for all $s \notin \Sigma$. That $W^s$ is closed follows directly from the a priori inequality of Proposition 2.3. Fix some elliptic Laplace–Beltrami operator $\Delta$ on $M$, and denote by $\Lambda^r$ the powers $\Lambda^r = (I - \Delta)^{r/2}$, which are elliptic pseudodifferential operators belonging to $S^r_{1,0}$. As a norm on $H^r$ we take $\|f\|_{H^r} = \|\Lambda^r f\|_{L^2}$. Since $W^s$ is closed in $H^s$, its codimension equals the dimension of

$$N_s = \{g \in H^{-s} : \langle g, Lu \rangle = 0 \text{ for all } u \in H^s \text{ such that } Lu \in H^s\},$$

where the pairing is the natural one of $H^{-s}$ with $H^s$. In particular, $g \in H^{-s}$ implies $\langle g, Lu \rangle = 0$ for all $u \in C^\infty(M)$, so that $L^*g = 0$ in the sense of distributions, where $L^*$ denotes the transpose of $L$. Thus the codimension of $W^s$ in $H^s$ is at most equal to the dimension of the nullspace $\tilde{N}_{-s}$ of $L^* = L$ in $H^{-s}$, in the sense of distributions.

We claim that if $0 \le s < s_0$ or $s_j < s < s_{j+1}$, then $\tilde{N}_{-s} \subset H^r$ for some $r > -s$. Since the inclusion of $H^r$ into $H^{-s}$ is a compact operator, this implies that $\tilde{N}_{-s}$ has finite dimension. The analysis of Section 6 and the first part of this section applies equally well to Sobolev spaces of negative order, so it suffices to verify that if $s \notin \Sigma$, then $-s$ is not in

$$\tilde{\Sigma} = \left\{s \in \mathbb{R} : \text{there exists } \tau \in \mathbb{R} \text{ such that } (s - \tfrac{1}{2} + i\tau)(s + \tfrac{1}{2} + i\tau) \in \Sigma_0\right\};$$

this is the definition of $\Sigma$ with the clause $s \in \mathbb{R}$ substituted for $s \in [0, \infty)$. But the substitution $(s, \tau) \to (-s, -\tau)$ leaves the quantity $(s - \tfrac{1}{2} + i\tau)(s + \tfrac{1}{2} + i\tau)$ invariant, so $s \in \Sigma$ if and only if $-s \in \Sigma$. $\square$

## 8. Remarks on Diffusions

In this section we reinterpret the global $C^\infty$ irregularity of $L$ in terms of the semigroup $e^{-\tau L}$, that is, in terms of the initial value problem

(8.1)
$$\begin{cases} (\partial_\tau + L)u(x,\tau) = 0, \\ u(x,0) = f(x) \end{cases}$$

for $(x,\tau) \in M \times [0,\infty)$. We assume that $L$ is of second order with $C^\infty$ coefficients and nonnegative principal symbol, elliptic on $M \setminus J$, equal to its formal adjoint in $L^2(M,\sigma)$, and of the form $-\partial_x^2 - \partial_t a^2(x,t)\partial_t$ in $V_{\delta_0}$, where $a(x,t)$ is as before and $d\sigma \equiv dx\, dt$ in $V_{\delta_0}$. We assume further that $L$ annihilates constant functions.

Denote by $L^2 \ominus 1 = \{u : u \perp 1\}$ the orthocomplement of the constant functions. The proofs of Lemmas 8.1 and 8.2 will be outlined at the end of the section. All norms without subscript denote the $L^2(M)$ norm, all functions are assumed real-valued, and $\langle f, g\rangle = \int_M fg\, d\sigma$.

**Lemma 8.1.** *There exists $C < \infty$ such that*
$$\|u\|^2 \leq C \langle Lu, u\rangle \quad \text{for all } u \in C^2 \text{ such that } u \perp 1.$$

As in Section 4, it follows that $L$ is essentially selfadjoint and maps the orthocomplement of 1 in its domain bijectively to $L^2 \ominus 1$. $\{e^{-\tau L} : \tau \geq 0\}$ is then a well-defined semigroup of contractions on $L^2(M)$, and there exists $c > 0$ such that $\|e^{-\tau L} f\| \leq e^{-c\tau} \|f\|$ for all $f \in L^2$ satisfying $f \perp 1$. Moreover, defining $L^{-1}$ to be the bounded linear operator on $L^2 \ominus 1$, with range in $\mathrm{Domain}(L) \ominus 1$, satisfying $LL^{-1} f = f$ for all $f \in L^2 \ominus 1$, we have
$$L^{-1} f = \int_0^\infty e^{-\tau L} f\, d\tau \quad \text{for all } f \in L^2 \ominus 1.$$

**Lemma 8.2.** *For each $s \in [0,\infty)$, each operator $e^{-\tau L}$ maps $H^s$ boundedly to $H^s$, and there exists $C < \infty$, depending only on $s$, such that*
$$\|e^{-\tau L}\|_{H^s} \leq Ce^{C\tau} \|f\|_{H^s} \quad \text{for all } f \in H^s.$$

This is a general fact that does not depend on the particular nature of $L$. For elliptic or subelliptic operators $L$ there is of course for each $s$ a far stronger inequality $\|e^{-\tau L} f\|_{H^s} \leq Ce^{-c\tau} \|f\|_{L^2}$ for all $f \in L^2 \ominus 1$.

**Proposition 8.3.** *There exist $f \in C^\infty(M)$, $s \in \mathbb{R}^+$, $\delta > 0$, and a sequence $\tau_j \to \infty$ such that for all $j$,*
$$\|e^{-\tau_j L} f\|_{H^s} \geq e^{+\delta \tau_j}.$$

If this were not so, $L$ would be globally regular in $C^\infty$. Indeed, consider any $f \in C^\infty$ satisfying $f \perp 1$. Consider any $r \in \mathbb{R}^+$, and fix $s > r$. Then, by hypothesis, $\|e^{-\tau L}f\|_{H^s} \leq C_\varepsilon e^{\varepsilon\tau}$ for all $\tau \geq 0$ and all $\varepsilon > 0$. Interpolating with the bound $Ce^{-c\tau}$ for the $L^2$ norm gives $\|e^{-\tau L}f\|_{H^r} \leq Ce^{-\delta\tau}$ for some $C, \delta \in \mathbb{R}^+$. Therefore $L^{-1}f = \int_0^\infty e^{-\tau L}f\,d\tau$ belongs to $H^r$. Since $r$ was arbitrary, we have $L^{-1}f \in C^\infty$, so $L$ would be globally regular. But this is not so, by the same reasoning as for the operators $L$ considered in the preceding sections.

*Sketch of proof of Lemma 8.1.* Integration by parts gives

$$\langle Lu, u\rangle \sim \|\nabla u\|^2_{L^2(M\setminus V_{\delta_0/4})} + \int_{V_{\delta_0}} |\partial_x u|^2 + |a\partial_t u|^2,$$

uniformly for all $u \in C^2(M)$, in the sense that for nonconstant $u$, the ratios of the two sides are bounded above and below by positive finite constants. Fix any closed ball $B \subset M$ disjoint from $V_{\delta_0}$, of positive radius. Define $N = \{u \in C^2(M) : \langle u, \varphi\rangle = 0\}$, where $\varphi \in C_0^2(B)$ satisfies $\int \varphi\,d\sigma \neq 0$. Then

$$\|u\|_{L^2(B)} \leq C\|\nabla u\|_{L^2(B)} \quad \text{for all } u \in N;$$

this variant of the usual Poincaré inequality is proved in the same way as the standard version. Clearly

$$\|u\|_{L^2(M\setminus V_{\delta_0/2})} \leq C\|\nabla u\|_{L^2(M\setminus V_{\delta_0/4})} + C\|u\|_{L^2(B)},$$

since $M\setminus V_{\delta_0/4}$ is connected. Thus

$$\|u\|^2_{L^2(M\setminus V_{\delta_0/2})} \leq C\langle Lu, u\rangle \quad \text{for all } u \in N.$$

But again it is clear that

$$\|u\|_{L^2(V_{\delta_0/2})} \leq C\|\partial_x u\|_{L^2(V_{\delta_0})} + C\|u\|_{L^2(V_{\delta_0}\setminus V_{\delta_0/2})},$$

so

$$\|u\|^2 \leq C\langle Lu, u\rangle \quad \text{for all } u \in N.$$

Given $f \in C^2$, we may decompose $f = u + \gamma$ in a unique way, with $u \in N$ and $\gamma \in \mathbb{R}$. Then $Lu = Lf \perp 1$, so

$$\|u\|^2 \leq C\langle Lu, u\rangle = C\langle Lf, u\rangle = C\langle Lf, f\rangle.$$

If $f \perp 1$, then $\gamma = -\sigma(M)^{-1}\int_M u\,d\sigma$, so $|\gamma| \leq C\|u\|$. Thus $\|f\| \leq \|u\| + |\gamma| \leq C\langle Lf, f\rangle$. □

*Sketch of proof of Lemma 8.2.* We begin with the derivation of an a priori upper bound for the norm of $e^{-\tau L}$ as an operator on $H^s$, assuming that $\tau \to e^{-\tau L}$ is a strongly continuous semigroup of bounded operators on $H^s$. Recall the formal identity

$$(8.2) \qquad [D, e^{-\tau L}] = -\int_0^\tau e^{-(\tau-t)L}[D, L]e^{-tL}\, dt,$$

for any linear operator $D$, which may be obtained by differentiating with respect to $\tau$ and writing $L[L, D] = [L, LD]$ to simplify the right-hand side.

It suffices, by interpolation, to treat the case where $s$ is a positive integer. Fix a basis $\{D_j\}$ for the $C^\infty(M)$ module of all differential operators of order at most $s$ on $M$, with $C^\infty$ coefficients. Write $\|T\|_{r,s}$ for the operator norm of $T : H^r \to H^s$. It suffices for our purpose to bound $\sum_j \|D_j \circ e^{-\tau L}\|_{s,0}$. But $e^{-\tau L} \circ D_j$ is bounded from $H^s$ to $H^0$, uniformly in $\tau$, so it suffices to bound

$$\sum_j \|[D_j, e^{-\tau L}]\|_{s,0} \leq \sum_j \int_0^\tau \|e^{-(\tau-t)L}[D_j, L]e^{-tL}\|_{s,0}\, dt.$$

Each $[D_j, L]$ may be expressed as a finite sum $\sum_i X_{ij}D_i$, where $X_{ij}$ is either a $C^\infty$ vector field that belongs to the span of $\{\partial_x, a\partial_t\}$ in $V_{\delta_0}$, or is multiplication by a $C^\infty$ function. Moreover, for $0 < \tau \leq 1$,

$$\|e^{-\tau L}X\|_{0,0} \leq C\tau^{-1/2}$$

for any such vector field $X$, since $\|Xf\|^2 \leq C\langle Lf, f\rangle$ for all $f$ and $\|e^{-\tau L}L^{1/2}\| = O(\tau^{-1/2})$. Thus, for $0 < \tau \leq 1$, we have

$$\sum_j \|[D_j, e^{-\tau L}]\|_{s,0} \leq C\sum_i \int_0^\tau |\tau - t|^{-1/2}\|D_i e^{-tL}\|_{s,0}\, dt$$

$$\leq C\tau^{1/2} + C\int_0^\tau |\tau - t|^{-1/2} \sum_i \|[D_i, e^{-tL}]\|_{s,0}\, dt.$$

Setting

$$z_T = \sup_{0 < \tau \leq T} \sum_j \|[D_j, e^{-\tau L}]\|_{s,0},$$

we thus have

$$z_T \leq CT^{1/2} + CT^{1/2}z_T,$$

where $C$ depends on $s$ and of course on $L$. Choosing $T$ to be sufficiently small therefore yields an a priori upper bound on $\|e^{-\tau L}\|_{s,s}$ for $0 < \tau \leq T$. The semigroup property then immediately implies $\|e^{-\tau L}\|_{s,s} \leq Ce^{C\tau}$ for all $\tau > 0$, for some $C < \infty$.

To make this formal argument rigorous, replace $L$ by $L^{(\varepsilon)} = L + \varepsilon\Delta$, where $\Delta$ is some elliptic selfadjoint second-order operator on $M$ with positive principal symbol, annihilating constants. The semigroup $\exp(-\tau L^{(\varepsilon)})$ is for each $\varepsilon > 0$ a strongly continuous semigroup of uniformly (in $\tau$, for fixed $\varepsilon$) bounded operators on each Sobolev space $H^s$. The above reasoning then applies for each $\varepsilon > 0$, uniformly in $\varepsilon$, using the observation that

$$\|e^{-\tau L^{(\varepsilon)}} \circ \varepsilon\nabla\|_{0,0} \leq C\varepsilon^{1/2}\tau^{-1/2} \leq C\tau^{-1/2}$$

because

$$\|\varepsilon^{1/2}\nabla u\|^2 \leq C\varepsilon\langle\Delta u, u\rangle \leq C\langle L^{(\varepsilon)} u, u\rangle.$$

The conclusion is that $\exp(-\tau L^{(\varepsilon)})$ is bounded on $H^s$, with a uniform bound as $\varepsilon \to 0$, for all $0 \leq \tau \leq 1$. An appropriate weak limit argument recovers the $H^s$ boundedness of $e^{-\tau L}$ from the uniform boundedness of the semigroup $\exp(-\tau L^{(\varepsilon)})$. $\square$

## 9. Concluding Remarks

1. The assumption that $a(x,t)$ vanishes to order precisely one as $t \to 0$, for $|x| \leq 1$, is essential. If $L$ satisfies hypotheses (2.1), (2.3) and (2.4) but $\partial a/\partial t(x,0) \equiv 0$ for all $|x| \leq 1$, then $L$ is globally regular in $C^\infty$ and moreover is exactly regular. This follows from the method of Boas and Straube [BS]. Indeed, the vector field $T = \partial_t$ satisfies $[T, \partial_x] \equiv 0$, while $[T, a\partial_t] = (\partial_t a) \cdot \partial_t$ vanishes identically on the set $J$ of all points where $L$ is not elliptic. The existence of a vector field $T$ for which both these commutators belong to the span of $\{\partial_x, a\partial_t\}$ on $J$, and such that $\{T, \partial_x, a\partial_t\}$ span the tangent space to $M$ at each point of $J$, suffices to imply global regularity.

2. From our point of view this holds because $\Lambda^s L \Lambda^{-s} - L$ takes the form $A_1 \circ \partial_x + A_2 \circ a\partial_t + A_3$, where $A_j \in S^0_{1,0}$ and $\sigma_0(A_j) \equiv 0$ on $\Gamma$. Thus $\Lambda^s L \Lambda^{-s}$ becomes an arbitrarily small perturbation of $L$ modulo relatively compact terms, for every $s$.

3. Suppose that (2.2) is replaced for some $m > 1$ by

$$(9.1) \qquad \frac{\partial^m a}{\partial t^m}(x,0) \neq 0 \quad \text{for all } |x| \leq 1,$$

but $\partial_t^j a(x,0) = 0$ for $|x| \leq 1$, for all $j < m$. Separation of variables, for the Dirichlet problem on $[-1,1] \times \mathbb{R}$ for the associated limiting operator $-\partial_x^2 - \alpha^2(x)(t^m\partial_t)^2 + \beta$, leads to solutions

$$(9.2) \qquad h(x)\exp(-\lambda t^{1-m})\chi_{t>0}$$

where $h \in C^\infty$ and $\lambda$ has positive real part. These solutions belong to $C^\infty$ for all $m > 1$, as is consistent with the preceding remarks.

4. Higher order vanishing of $a(x,t)$ at $t = 0$ may nonetheless be related to more subtle singularities of solutions of $L$. Suppose that $M$ is a real analytic manifold, that (9.1) holds for some integer $m > 1$ and that $a$ belongs to every Gevrey class $G^s$ of order $s$ strictly greater than one.

**Question.** Is it true, perhaps under some subsidiary hypotheses on $a$, that if $m > 1$ then $L$ is globally regular in $G^r$ if and only if $r \geq m/(m-1)$?

This is suggested by two observations. Firstly, the coefficient $h$ in (9.2) belongs to $G^s$ for every $s > 1$, and hence if $h$ does not vanish identically, the solution defined by (9.2) belongs to $G^r$ if and only if $r \geq m/(m-1)$.

Secondly, if $\Lambda = \exp(|\partial_t|)^p$, then formally
$$[\Lambda, t^m \partial_t]\Lambda^{-1} = ct^{m-1}|\partial_t|^p = c(t^m|\partial_t|)^{(m-1)/m}$$
if $p = (m-1)/m$. Thus $\Lambda L \Lambda^{-1} - L$ differs from $L$ by various terms, of which the worst is formally $O(t^m|\partial_t|)^{2-m^{-1}}$, which should be bounded relative to $(t^m \partial_t)^2$. An $L^2$ estimate for $\Lambda L \Lambda^{-1}$ would imply global $G^r$ regularity for $L$, with $r = m/(m-1)$.

5. Exact regularity in the scale of Sobolev spaces has played an important part in our analysis. The following question has been asked by J. J. Kohn in the context of the $\bar{\partial}$–Neumann problem.

**Question.** Does global regularity in $C^\infty$ always imply exact regularity, under the hypothesis $\|u\| \leq C\|Lu\|$?

The answer is certainly "no" for general operators, if the basic $L^2$ estimate $\|u\| \leq C\|Lu\|$ is not assumed. For instance, the operators $\partial_{x_1} + \alpha \partial_{x_2}$ mentioned in the introduction are counterexamples, as are elliptic pseudodifferential operators of negative order. On the other hand, no counterexamples are presently known for the $\bar{\partial}$–Neumann problem, where the $L^2$ estimate holds automatically. The following more narrowly focused question is a natural generalization of the $\bar{\partial}$–Neumann situation.

**Question.** Assume $L$ to be selfadjoint and of second order, with nonnegative principal symbol and $C^\infty$ coefficients. Suppose that $\|u\| \leq C\|Lu\|$ for all $u$ in the domain of $L$. Does global regularity for $L$ in $C^\infty$ imply exact regularity? If not, does this hold under some natural and weak additional hypotheses?

Michael Christ, University of California, Los Angeles, Los Angeles, CA 90095-1555

*E-mail address*: `christ@@math.ucla.edu`

*Current address*: Mathematical Sciences Research Institute, 1000 Centennial Road, Berkeley, CA 94720-5070